\documentclass{article}
\usepackage{amssymb}
\usepackage{amsmath}
\usepackage{amsfonts}

\newtheorem{theorem}{Theorem}

\newtheorem{definition}[theorem]{Definition}

\input{tcilatex}

\begin{document}

\title{Samuelson's webs of maximum rank}
\author{Vladislav V. Goldberg and Valentin V. Lychagin}
\date{}
\maketitle

\begin{abstract}
The authors found necessary and sufficient conditions
for Samuelson's web to be of maximum rank.
\end{abstract}


\section{Introduction and Motivations}

A planar 3-web $W$ can be defined by three differential $1$-forms, say, $%
\omega_{1}, \omega_{2}$ and $\omega_{3}$, where $\omega_{1}\wedge \omega_{2}\neq
0, \;\; \omega_{2}\wedge \omega_{3}\neq 0$
and $\omega_{1}\wedge
\omega_{3}\neq 0.$ These forms can be normalized in such a way that $%
\omega_{1} + \omega_{2} +\omega_{3}=0$. They satisfy the following structure
equations:
\begin{equation*}
d \omega_{1} = \omega_{1} \wedge \gamma, \;\; d\omega_{2}= \omega_{2} \wedge
\gamma, \;\; d \gamma = K \omega_{1} \wedge \omega_2,
\end{equation*}
where $\mathcal{K} (W) = d \gamma$ is the \emph{web curvature $2$-form}, and $K$ is
its \emph{scalar curvature}. If $f (x, y)$ is a web function, then (see \cite{B55})
\begin{equation*}
K = - \frac{1}{f_x f_y}
\Biggl(\log \biggl(\frac{f_x}{f_y}\biggr) \Biggr)_{xy}.
\end{equation*}

The condition $\mathcal{K} (W) = 0$ is necessary and sufficient for a 3-web
to be parallelizable (trivial), i.e., to be equivalent to a 3-web formed by
three foliations of parallel lines of an affine plane $\mathbb{A}^2$.

This part of the web theory was used by Gerard Debreu, a Nobel Prize winner in Economics (1983). In \cite{D54} Debreu obtained conditions for a preference ordering to be representable by a numerical function. After
proving that such a function exists, Debreu in \cite{D59} and \cite{D60} investigated when it would be additively separable and proved that this question is equivalent to requiring that a planar 3-web given by the level curves
of the function, the verticals and the horizontals be
equivalent to the trivial 3-web (see a more extensive
treatment in \cite{V02} and \cite{W89}).

This required the satisfaction of the \emph{hexagon condition} \cite{B55}.
If $\mathcal{K} (W) \neq 0$, then there is an obstruction to triviality of a planar 3-web, i.e., the failure of a hexagon consisting of "threads" of the web to be closed. Russell \cite{Ru03} showed
that $\mathcal{K} (W)$ measures the local failure of
the hexagon to close (in economic terms,
$\mathcal{K}(W)$ is the local failure of the Expected
Utility Maximization (EUM) Axioms). As was indicated in
\cite{Ru97}, this 2-form $\mathcal{K} (W)$ was first identified by Pareto in \cite{P71}. The EUM hypothesis with limited experimental data was tested in \cite{CRS04}.

Note that Samuelson \cite{S02} derived a third-order PDE which is equivalent to the condition $\mathcal{K} (W) = 0.$ Economists called this PDE Samuelson's equation. However, in fact this equation is well-known St. Robert equation (see, for example, \cite{AS92}, p. 43).

Another application of the web theory in economics was given by Paul A. Samuelson, a Nobel Prize winner in Economics (1970). Using the web theory terms, we can say that Samuelson asked for an analytic criterium for a
certain 2-web to satisfy a certain natural area condition.

In \cite{CR08} the authors present an overview of this application of the web theory to economics.

If a 2-web is formed by the level sets of two functons $u (x, y) = %
\mbox{const.}, v (x, y) = \mbox{const.}$,
 then the area condition is
\begin{equation}
\frac{a}{c}=\frac{b}{d},  \label{area}
\end{equation}
where $a, b, c,$ and $d$ are the areas of quadrilaterals bounded by "threads" of the web. We shall call (\ref{area}) \emph{Samuelson's area condition} (or $S$-\emph{condition}).

Note that area condition (\ref{area}) has been used
by J. C. Maxwell (see \cite{M98}) in his classic work on thermodynamics.

The area condition was studied in detail in \cite{S60},
\cite{S72}, \cite{S83}, \cite{CRS01}, and \cite{CR06}.

Hess \cite{H80} considered Lagrangian 2-webs under the name of bipolarized symplectic manifolds and introduced a connection which in the planar case measures the failure of the equality $ab = cd$.

For a more detailed discussion of the properties of Hess'
connection see \cite{Va85} and \cite{Va89}.

A bi-Lagrangian manifold is a symplectic manifold endowed with two natural Lagrangian foliations. It was recently investigated in detail in \cite{EST}.
This manifold admits Hess' connection \cite{H80}.

In his Nobel lecture \cite{S72} Samuelson used
the area condition (\ref{area}) to characterize \emph{profit maximization}. His test for profit maximization is
as follows: calculate the areas $a, b, c,$ and $d$ of any four quadrilaterals cut out by the leaves of the demand systems. Then if $ad = bc$, then the firm is maximizing profits.

Tabachnikov has shown in \cite{Ta93} that when the Samuelson area condition is satisfied, the $2$-web is trivial under an area preserving transformation. Thus we can calibrate the leaves of the web in such a way that there is unit area between the threads labeled $x$ and $x+1$ and $y$ and $y+1$, for each respective member of the family. In classical thermodynamics this calibration corresponds to the passage from empirical temperature and entropy to absolute temperature and entropy. In economics this recalibration is not possible, so the relevant test for maximization is whether or not the already calibrated $2$-web, when trivialized to the horizontal/vertical web, already satisfies the equal area condition.

Tabachnikov \cite{Ta93} gives a further characterization of the profit maximizing condition. If we place the standard area (symplectic) form $\mbox{\bf d}\mbox{\bf L}_1\wedge\mbox{\bf d} \mbox{\bf W}_1$ on $\mathbb{R}^2$, the $2$-web of factor demands is a Lagrangian $2$-web since all curves on the symplectic plane are Lagrangian
submanifolds.

 Theorem 0.1 in \cite{Ta93} now applies directly and a symplectic, torsion-free connection (the Hess connection \cite{H80}) can be associated with the web. As Tabachnikov shows, when the Samuelson test is satisfied with equality, this connection is flat. This provides an alternative characterization of profit maximizing behavior.

Finally note that if one imposes the standard symplectic form
$$ \mbox{\bf dL}_1\wedge \mbox{\bf dW}_1+\mbox{\bf dL}_2\wedge \mbox{\bf   dW}_2
$$
on $\mathbb{R}^4$, the $2$-dimensional submanifold  is a Lagrangian submanifold of $\mathbb{R}^4$. This means that the mapping from $\mbox{\bf L}_1, \mbox{\bf W}_1$ to $\mbox{\bf L}_2, \mbox{\bf W}_2$ is area preserving and orientation reversing. Since maximizing economic processes take place on Lagrangian submanifolds, economics, too, succumbs to Weinstein's Lagrangian creed. That everything is a Lagrangian submanifold \cite{We81}.

\section{Samuelson's Webs}

Let $M$ be a two-dimensional manifold, and let a planar $4$-web $<\omega_1, \omega_2,\omega_3,\omega_4>$ be given  on $M$.

\begin{definition} Such a $4$-web is said to be the \textbf{Samuelson's web}, if the forms
$\omega_i, i= 1, 2,  3, 4,$ satisfy the following exterior quadratic relation:
\begin{equation}
\omega_{3}\wedge \omega_{1}+\omega_{4}\wedge \omega_{2}=0.
\label{S-condition}
\end{equation}
\end{definition}
In what follows, for brevity we shall call Samuelson's webs $S$\emph{-webs}.

Consider now the main example of an $S$-web.

Let $\mathbb{R}^4$ be a four-dimensional symplectic manifold with a structure form $dy_1 \wedge dx_1 + dy_2 \wedge dx_2$, and let
$M^2 \subset \mathbb{R}^4$ be such a Lagrangian surface on which any pair of the coordinate functions $x_1, x_2, y_1, y_2, x_i, y_j$ is functionally independent. Then the $4$-web on this surface defined by the level curves of these functions is an $S$-web.

Let $M$ be an $2$-manifold. In Definition 1 the
forms $\omega_i$ are defined up to factors $\lambda_i$:
\begin{equation}
\omega_i \rightarrow \lambda_i \omega_i,
\label{lambdai}
\end{equation}
 which satisfy
the following condition:
\begin{equation}
\lambda_{3}\lambda_{1}=\lambda_{4}\lambda_{2}.
\label{lambda-condition}
\end{equation}

Observation (\ref{lambda-condition}) allows us to make
the following normalization of a 4-web: we can take
the factors $\lambda_{1}, \lambda_{2}$ and $\lambda_{3}$
in such a way that
\begin{equation}
\omega_{3} + \omega_{1} + \omega_{2}=0.
\label{1st normalization}
\end{equation}
Under this choice of $\omega_{1}, \omega_{2}$ and
$\omega_{3}$, taking into account (\ref{S-condition}),
we can prove that the factors $\lambda_i$ in $\omega_i$ in
(\ref{lambdai}) must be equal:
\begin{equation}
\lambda_{1} = \lambda_{2} = \lambda_{3}=\lambda.
\label{lambda}
\end{equation}
By (\ref{1st normalization}), $S$-condition (\ref{S-condition})
becomes
$$
(\omega_{4} + \omega_{1}) \wedge \omega_{2} = 0.
$$
It follows that
\begin{equation}
\omega_{4} + \omega_{1} + b\omega_{2}=0.
\label{2nd normalization}
\end{equation}

We shall call the normalization
\begin{equation}
\renewcommand{\arraystretch}{1.5}
\begin{array}{ll}
\omega_{1} + \omega_{2} + \omega_{3}=0, \\
\omega_{4} + \omega_{1} + b\omega_{2}=0
\end{array}%
\renewcommand{\arraystretch}{1}
\label{canonical normalization}
\end{equation}
\emph{canonical}.

Let us take
$$
\omega_3 = df, \;\; \omega_1 \wedge dx = 0, \;\; \omega_2 \wedge dy = 0.
$$
Then the functions $x$ and $y$ can be viewed as coordinates in the plane, and the first equation of
(\ref{canonical normalization}) gives
\begin{equation}
\omega_3 = df, \;\; \omega_1 = - f_x dx, \;\; \omega_2 = - f_y dy.
\label{omegai}
\end{equation}

Let $g$ be another function with $\omega_4 \wedge dg = 0$. Then $\omega_4 = \lambda dg$, and the second
equation of
(\ref{canonical normalization}) gives
$$
\lambda dg - f_x dx - f_y dy = 0.
$$
Therefore,
$$
\lambda g_x = f_x, \;\; \lambda g_y = b f_y,
$$
and
$$
b= \frac{f_x g_y}{f_y g_x},\;\;\lambda=\frac{f_x }{g_x}.
$$
As a result, we have
$$
\omega_4 = \lambda dg = f_x dx + \frac{f_x}{g_x} g_y dy,
$$
or
\begin{equation}
\omega_4 = f_x dx + b f_y dy. \label{omega4}
\end{equation}

It is easy to see that by (\ref{omegai}) and (\ref{omega4}) relation
(\ref{S-condition}) is satisfied:

$$
\renewcommand{\arraystretch}{1.5}
\begin{array}{ll}
\omega_{3} \wedge \omega_{1} + \omega_4 \wedge \omega_{2} = (f_x dx + f_y dy) \wedge (-f_x dx) + (f_x dx + b f_y dy) \wedge (-f_y dy) \\= f_x f_y dx \wedge dy - f_x f_y dx \wedge dy = 0.
\end{array}
\renewcommand{\arraystretch}{1}
$$

\section{Structure Equations}
As in \cite{GL06}, we denote by $\gamma$ such 1-form that $$
d\omega_i = \omega_i \wedge \gamma, \;\; i = 1, 2,  3.
$$
The form $\gamma$ defines the Chern connection in the plane. The curvature 2-form of this connection $d\gamma$ is an invariant of the
$3$-web $<\omega_1, \omega_2,\omega_3>$.

Moreover, we find that
\begin{equation}
\gamma = H \omega_3, \label{gamma}
\end{equation}
where
\begin{equation}
H = \frac{f_xy}{f_xf_y} \label{H}
\end{equation}
and
\begin{equation}
d\gamma = K \omega_1 \wedge \omega_2, \label{dgamma}
\end{equation}
where
\begin{equation}
K = - \frac{1}{f_x}{f_y} \biggl(log \bigl(\frac{f_x}{f_y}\bigr)\biggr)_{xy} \label{K}
\end{equation}
is the \emph{scalar curvature} of the connection (or a $3$-web).

Denote by $\partial_1$ and  $\partial_2$ the basis of vector fields which is dual to the cobasis $\{\omega_1, \omega_2\}$:
$$
<\omega_i, \partial_j> = \delta_{ij}, \;\; i. j = 1, 2.
$$
Then for any function $p$ one has
\begin{equation}
dp = p_1 \omega_1 + p_2 \omega_2, \label{dp}
\end{equation}
where $p_1 = \partial_1 (p)$ and $p_2 = \partial_2 (p)$.

Taking differential of (\ref{dp}), we get
$$
\renewcommand{\arraystretch}{1.5}
\begin{array}{ll}
0 = dp_1 \wedge \omega_1 + dp_2 \wedge \omega_2
+ p_1 \omega_1 \wedge \gamma  + p_2 \omega_2 \wedge \gamma \\
= - \partial_2 \partial_1 (p) \omega_2 \wedge \omega_1
+ \partial_1 \partial_2 (p) \omega_1 \wedge \omega_2 +
H (p_1 - p_2) \omega_1 \wedge \omega_2,
\end{array}%
\renewcommand{\arraystretch}{1}
$$
or
\begin{equation}
[\partial_1, \partial_2] = H (\partial_2 - \partial_1).
\label{commutator}
\end{equation}
Remark that
\begin{equation}
K = \partial_1 (H) - \partial_2 (H).
\label{Kbis}
\end{equation}

In the coordinates $(x, y)$, the vector fields
$\partial_1$ and $\partial_2$ have the following form:
$$
\partial_1 = - \frac{1}{f_x} \partial_x, \;\; \partial_2 = - \frac{1}{f_y} \partial_y,
$$
and (\ref{commutator}) can be verified by direct calculation.

In what follows we shall use the following notation:
$$
p_i = \partial_i (p), \;\; p_{ij} = \partial_i \partial_j (p), \;\; \textrm{etc}.
$$

\section{Samuelson's Equations}
Samuelson's condition (\ref{S-condition}) means that there are positive factors $s_1, s_2, t_1$ and $t_2$ such that the forms  $s_1 \omega_1, s_2 \omega_2, t_1\omega_3$ and $t_2 \omega_4$ satisfy (\ref{S-condition}) and are closed.

These conditions imply the following relations:
\begin{equation}
\renewcommand{\arraystretch}{1.5}
\left\{
\begin{array}{ll}
d(s_1 \omega_1) = d(s_2 \omega_2) = d(t_1\omega_3)
= d(t_2 \omega_4) = 0, \\
s_1 t_1 = s_2 t_2.
 \label{S-equation}
\end{array}%
\right.
\renewcommand{\arraystretch}{1}
\end{equation}

We derive now an explicit form of Samuelson's equations (\ref{S-equation}). We have
$$
d (s_1 \omega_1) = ds_1 \wedge \omega_1 + s_1 \omega_1 \wedge \gamma = -s_{1,2} \omega_1 \wedge \omega_2 + H s_1 \omega_1\wedge \omega_2 = 0.
$$
It follows that
$$
s_{1,2} = H s_1.
$$

Similarly, we have
$$
d (s_2 \omega_2) = ds_2 \wedge \omega_2 + s_2 \omega_2 \wedge \gamma = s_{2,1} \omega_1 \wedge \omega_2 - H s_2 \omega_1\wedge \omega_2 = 0,
$$
or
$$
s_{2,1} = H s_2.
$$

For the third equation of (\ref{S-equation}) one has
$$
d (t_1 \omega_3) = -dt_1 \wedge (\omega_1 + \omega_2) - t_1 (\omega_1 + \omega_2) \wedge \gamma = (t_{1,2} - t_{1,1}) \omega_1\wedge \omega_2 = 0,
$$
or
$$
t_{1,2} - t_{1,1} = 0.
$$
Similarly, we have
$$
\renewcommand{\arraystretch}{1.5}
\begin{array}{ll}
d (t_2 \omega_4) = -dt_2 \wedge (\omega_1 + b \omega_2) - t_2 (d\omega_1 + bd\omega_2 + db \wedge \omega_2) \\
= t_{2,2} \omega_1\wedge \omega_2 -b t_{2,2}) \omega_1\wedge \omega_2 - t_2 H \omega_1\wedge \omega_2 + t_2 b H \omega_1\wedge \omega_2 - t_2 b_1 \omega_1\wedge \omega_2= 0,
\end{array}%
\renewcommand{\arraystretch}{1}
$$
or
$$
t_{2,2} - bt_{2,1} - t_2(b_1 - (b-1)H) = 0.
$$

Define the new functions $\sigma_i$ and $\tau_i, \, i = 1, 2$:
\begin{equation}
s_i = \log |\sigma_i|, \;\; t_i = \log |\tau_i|. \label{sigma/tau}
\end{equation}
Then the Samuelson equations take the following form:
\begin{equation}
\renewcommand{\arraystretch}{1.5}
\left\{
\begin{array}{ll}
\sigma_{1,2} = H, \;\; \sigma_{2,1} = H, \\
\tau_{1,2} - \tau_{1,1} = 0, & \\
 b \tau_{2,1} - \tau_{2,2} = (b - 1)H - b_1,
 \label{S-equations}
\end{array}%
\right.
\renewcommand{\arraystretch}{1}
\end{equation}

In addition, the second equation of (\ref{S-equation}) and (\ref{sigma/tau}) imply that
\begin{equation}
\sigma_1 + \tau_1 = \sigma_2 + \tau_2. \label{sigma+tau}
\end{equation}

Taking into account the last equation of (\ref{S-equations}), representing $\tau_2$ from
 (\ref{sigma+tau}) in the form $\tau_2 = \sigma_1 + \tau_1 - \sigma_2$ and applying (\ref{commutator}), we get the final form of
 Samuelson's equations:
\begin{equation}
\renewcommand{\arraystretch}{1.5}
\left\{
\begin{array}{ll}
\sigma_{1,2} = H, \;\; \sigma_{2,1} = H, \\
\tau_{1,2} - \tau_{1,1} = 0, & \\
 b \sigma_{1,1} + (b-1) \tau_{1,2}+ \sigma_{2,2} = 2bH - b_1.
 \label{final S-equations}
\end{array}%
\right.
\renewcommand{\arraystretch}{1}
\end{equation}

We compute now the first and second prolongations of PDE system (\ref{final S-equations}). For the first equation of (\ref{final S-equations}) we have
\begin{equation}
\renewcommand{\arraystretch}{1.5}
\left\{
\begin{array}{lll}
\sigma_{1,2} = H, && \\
\sigma_{1,12} = H_1, & \sigma_{1,22} = H_2, & \\
\sigma_{1,112} = H_{11}, & \sigma_{1,122} = H_{12}, &
\sigma_{1,222} = H_{22}.
 \label{prolongation of 1st eq}
\end{array}%
\right.
\renewcommand{\arraystretch}{1}
\end{equation}

 For the second equation of (\ref{final S-equations}) we have
\begin{equation}
\renewcommand{\arraystretch}{1.5}
\left\{
\begin{array}{lll}
\sigma_{2,1} = H, && \\
\sigma_{2,11} = H_1, \;\; \sigma_{2,12} = H_2 - H^2 + H \sigma_{2,2}, && \\
\sigma_{2,111} = H_{11}, \;\; \sigma_{2,112} = H_{12} - 2 HH_1 +HH_2 - H^3 +(H^2 +H_1)\sigma_{2,2}, &&\\
\sigma_{2,122} = 2H \sigma_{2,22} + (H_2 - H^2) \sigma_{2,2} + H_{2,2} - 3HH_2 + H^3. &&
 \label{prolongation of 2nd eq}
\end{array}%
\right.
\renewcommand{\arraystretch}{1}
\end{equation}

Solving the third equation of (\ref{S-equation}),
we find that
$$
\tau_1 = w (f)
$$
for some function $w$.

Remark that by (\ref{omegai}), $\partial_i = - 1, \; i = 1, 2$, and as a result, $\partial_i (w) = - w'$.

Let us rewrite system (\ref{S-equation}) in the form
\begin{equation}
\renewcommand{\arraystretch}{1.5}
\left\{
\begin{array}{ll}
\sigma_{1,2} = H, \;\; \sigma_{2,1} = H, \\
\tau_{1,2} - \tau_{1,1} = 0, & \\
 b \sigma_{1,1} + \sigma_{2,2} = B,
 \end{array}%
\right.
\renewcommand{\arraystretch}{1}
\label{S-red}
\end{equation}
where
$$
B = 2bH - b_1 + (b-1) w'.
$$

We shall investigate the solvability
of (\ref{S-red}) with respect to $\sigma_1$ and
$\sigma_2$.

System (\ref{S-red}) is the first-order system of PDE, and its first prolongation has the form
\begin{equation}
\renewcommand{\arraystretch}{1.5}
\left\{
\begin{array}{ll}
\sigma_{1,12} = H_1, \;\; \sigma_{1,22} = H_2, \\
\sigma_{2,11} = H_1, \;\; \sigma_{2,21} = H_2, \\
\sigma_{2,11} = R + r \sigma_{1,1},\\
 \sigma_{2,22} = B_2 - bH_1 +bH^2-(b_1 + bH)\sigma_{1,1},
 \label{1st prol}
\end{array}%
\right.
\renewcommand{\arraystretch}{1}
\end{equation}
where
$$
R=\frac{b_1 - H_2 + H^2 - HB}{b}, \;\; r = H - \frac{b_1}{b}.
$$

Note that by (\ref{commutator}) the first and the third  equations of (\ref{1st prol}) imply that
\begin{equation}
\sigma_{1,21} = H_1 - H^2 + H \sigma_{1,1}, \;\;
\sigma_{1,12} = H_2 - H^2 + HB - Hb \sigma_{1,1}.
\end{equation}

Computing the third derivatives of $\tau_1$ and $\tau_2$, we shall get four relations if we use different ways of finding $\sigma_{1,112}$ and $\sigma_{1,122}$ as well as $\sigma_{2,112}$ and $\sigma_{2,122}$.

Denote  PDE system (\ref{S-red}) by $\mathfrak{E}_1 \in \mathfrak{J}^1 (\pi)$, where $\pi: \mathbb{R}^2: \mathbb{R}^2 \times \mathbb{R}^2$ is the trivial bundle, and denote the first prolongation (\ref{1st prol}) of (\ref{S-red}) by $\mathfrak{E}_1^{(1)} \in \mathfrak{J}^2 (\pi)$. Then we get the following tower:
\begin{equation}
\mathbb{R}^2 \overset{\pi}{\leftarrow} \mathfrak{J}^1\overset{\pi_{1,0}}\leftarrow \mathfrak{E}_1 \overset{\pi_{2,1}}\leftarrow\mathfrak{E}_1^{(1)},
\label{tower}
\end{equation}
where the map $\pi_{2,1}$ is the diffeomorphism,
and $\pi_{1,0}$ is the one-dimensional bundle with (as we saw) fiberwise coordinate $\sigma_{1,1}$. Therefore, as it follows from \cite{L88} and \cite{KL08}, we have the only obstruction for integrability, and this function can be found by different computations of the third derivatives.

In our case two different computations of $\sigma_{1,112}$ give us
$$
\sigma_{1,112} = \partial_1 (\sigma_{1,12}) = H_{1,1}
$$
and
$$
\sigma_{1,211} = \partial_2 (\sigma_{1,11}) = R_2 + r_2
\sigma_{1,11} + r \sigma_{1,21} \\
= R_2 + rH_1 +rH^2 + (r_2 +rH_1) \sigma_{1,1}.
$$
On the other hand, we have
$$
\renewcommand{\arraystretch}{1.5}
\begin{array}{ll}
\sigma_{1,211} &= \sigma_{1,121} + H(\sigma_{1,11} - \sigma_{1,21}) \\
&= \sigma_{1,121} +(RH - HH_1 + H^3 + 2H-H^2)\sigma_{1,1} \\ &= \sigma_{1,112} + \partial_1(H(\sigma_{1,11}-\sigma_{1,12}))+ (RH - HH_1 + H^3 + 2H-H^2)\sigma_{1,1}\\
&= H_{11} + H^3 - 3HH_1 + 2HR + (H_1 - H^2 + 2H_2)
\sigma_{1,1}.
\end{array}%
\renewcommand{\arraystretch}{1}
$$
Then the obstruction to the formal integrability
of (\ref{S-red}) vanishes if and only if the following two conditions are satisfied:
\begin{equation}
R_2 +rH_1 - rH^2 - H_{11} +3HH_1 - H^3 - 2HR = 0
\label{1st int cond}
\end{equation}
and
\begin{equation}
r_2 + H^2 - H_{1} -rH = 0.
\label{2nd int cond}
\end{equation}

Substituting $r = H - \frac{b_1}{b}$ into (\ref{2nd int cond}) and applying (\ref{Kbis}), we find that
\begin{equation}
K = - \delta_2 \delta_1 \log b.
\label{Kbisbis}
\end{equation}

\begin{definition}
The rank of the system of equations $(\ref{1st int cond})$ and $(\ref{2nd int cond})$ is said to be the \textbf{rank of Samuelson's web}.
\end{definition}

Remark that condition (\ref{2nd int cond}) does not contain the function $w$, and condition (\ref{1st int cond})
can be written in the form
\begin{equation}
T_3 w''' + T_2 w'' + T_1 w + T_0 = 0.
\label{1st int condbis}
\end{equation}
Given $w$, tower (\ref{tower}) shows that we can get a solution space of dimension not exceeding three. In order to get a three-dimensional solution space, we need both conditions  (\ref{1st int cond}) and (\ref{2nd int cond}). Condition (\ref{1st int condbis}) now is a third-order ODE with respect to $w$. We rewrite
(\ref{1st int condbis}) in the form
\begin{equation}
 w''' + \frac{T_2}{T_3} w'' + \frac{T_1}{T_3} w + \frac{T_0}{T_3} = 0.
\label{1st int condbisbis}
\end{equation}

If the coefficients of (\ref{1st int condbisbis}) depend on $f$ only, we have an extra three-dimensional solution space, and therefore the rank of the Samuelson web equals six. Otherwise, the rank of the Samuelson web is less than six.

Keeping in mind that our second condition (\ref{2nd int cond}) has the form (\ref{Kbisbis}), we have the following theorem:
\begin{theorem}
A Samuelson web has the maximum rank six if and only if the following conditions hold:
$$
\renewcommand{\arraystretch}{1.5}
\left\{
\begin{array}{ll}
K = - \delta_2\delta_1 \log b,\\
\delta\biggl({\displaystyle\frac{T_2}{T_3}}\biggr) =
\delta\biggl({\displaystyle\frac{T_1}{T_3}}\biggr) =
\delta\biggl({\displaystyle\frac{T_0}{T_3}}\biggr) = 0,
\end{array}
\right.
\renewcommand{\arraystretch}{1}
$$
where $\delta = \partial_1 - \partial_2$.
\end{theorem}

\end{document}